\documentstyle[11pt,leqno,amsfonts,amssymb]{article}
\newtheorem{definition}{Definition}[section]
\newtheorem{lemma}[definition]{Lemma}
\newtheorem{theorem}[definition]{Theorem}

\newtheorem{remark}[definition]{Remark}

 1

\newenvironment{proof}{\noindent{\bf Proof~:}}{\QED\medskip}
\def\QED{\hskip0.1em\hfill\null\ \null\nobreak\hfill
\kern3pt\lower1.8pt\vbox{\hrule\hbox
{\vrule\kern1pt\vbox{\kern1.7pt \hbox{$\scriptstyle
QED$}\kern0.2pt}\kern1pt\vrule}\hrule}}

\newcommand{\ta}{\tilde{a}}
\newcommand{\tb}{\tilde{b}}
\newcommand{\tc}{\tilde{c}}
\newcommand{\ha}{\hat{a}}
\newcommand{\hb}{\hat{b}}
\newcommand{\hc}{\hat{c}}

\newcommand{\ZZ}{{\mathbb{Z}}}

\newcommand{\LL}{{\mathbb{L}}}
\newcommand{\RR}{{\mathbb{R}}}

\newcommand{\iso}{\cong}

\newcommand{\inc}{\hookrightarrow}
\newcommand{\too}{\longrightarrow}

\newcommand{\la}{\langle}
\newcommand{\ra}{\rangle}
\newcommand{\bd}{\partial}

\begin{document}

\title{Non-formal compact manifolds with small Betti numbers}

\author{Marisa Fern\'andez
and Vicente Mu\~noz}

\date{May 10, 2006}

\maketitle

\begin{abstract}
We show that, for any $k\geq 1$, there exist non-formal compact
orientable $(k-1)$-connected $n$-manifolds with $k$-th Betti
number $b_k=b\geq 0$ if and only if $n\geq \max \{ 4k-1 ,
4k+3-2b\}$.
\end{abstract}

\section{Introduction} \label{introduction}
Simply connected compact manifolds of dimension less than or equal
to $6$ are formal \cite{NM}. Moreover $(k-1)$-connected compact
orientable
manifolds of dimension less than or equal to $4k-2$ are formal
\cite{Mi,FM}, for any $k\geq 1$. To check that this bound is
optimal, a method to construct non-formal simply connected compact
manifolds of any dimension $n\geq 7$ was given
in
\cite{FM2} (see also \cite{O} for a previous construction of an
example of dimension $n=7$). Later Dranishnikov and Rudyak
\cite{Ru} extended the result to give examples of $(k-1)$-connected
non-formal manifolds of dimension $n\geq 4k-1$.

For any $k\geq 1$, there is an alternative approach used by
Cavalcanti in \cite{Gil} to construct $(k-1)$-connected compact
orientable non-formal manifolds, which gives examples with small
Betti number $b_k$. If one focusses in the case $k=1$, in
\cite{FM3} it is studied how the smallness of $b_1$ may force
the formality of a manifold. More concretely, in \cite{FM3} it is
proved that any compact orientable $n$-dimensional manifold with
first Betti number $b_1=b$ is formal if $n\leq \max \{ 2, 6-2b
\}$, and that this bound is sharp, i.e., there are non-formal
examples of compact orientable manifolds if: \
  (a) $b_k=0$ and $n\geq 7$; \
  (b) $b_k=1$ and $n\geq 5$; \ or \
  (c) $b_k\geq 2$ and $n\geq 3$.
The examples constructed there follow the lines of those of
\cite{FM2}. For the general case $k\geq 1$, Cavalcanti \cite{Gil}
proves that $(k-1)$-connected compact orientable $n$-manifolds
with $b_k=1$ are formal whenever $n\leq 4k$.

The natural \textit{geography} question that arises in this
situation is the following:

\bigskip

\hspace{5mm}
\begin{minipage}{14cm} {For which $(n,k,b)$ with $n,k
\geq 1$, $b\geq 0$, are there compact orientable $(k-1)$-connected
manifolds of dimension $n$ and with $b_k=b$ which are non-formal?}
\end{minipage}

\bigskip

(Note that the orientability condition is only relevant if $k=1$.)
In this paper we solve completely the above problem by proving
the following main result.

\begin{theorem} \label{thm:main}
For any $k\geq 1$, there exist non-formal compact orientable
$(k-1)$-connected $n$-manifolds with $k$-th Betti number
$b_k=b\geq 0$ if and only if $n\geq \max \{ 4k-1 , 4k+3-2b\}$.
\end{theorem}

The above result can be restated as follows: let $M$ be a compact
orientable $(k-1)$-connected manifold of dimension $n$. Then $M$
is formal if:
 \begin{enumerate}
 \item[(a)] $b_k=0$ and $n\leq 4k+2$;
 \item[(b)] $b_k=1$ and $n\leq 4k$; or
 \item[(c)] $b_k\geq 2$ and $n\leq 4k-2$.
 \end{enumerate}
In all other situations, namely
 \begin{enumerate}
 \item[(a')] $b_k=0$ and $n\geq 4k+3$;
 \item[(b')] $b_k=1$ and $n\geq 4k+1$; or
 \item[(c')] $b_k\geq 2$ and $n\geq 4k-1$.
 \end{enumerate}
there are non-formal examples.

The examples that we construct in this paper follow the lines of
\cite{Ru} (see also Example 5 in \cite{FM} where the same
construction is used). Some (alternative) examples of compact
orientable $(k-1)$-connected $n$-manifold in the cases (b') and
(c') are given in \cite{Gil}, but the list is not exhaustive (for
instance it does not cover the case $b_k=1$ and $n=4k+1$, see
remark \ref{rem:Gil}).

\section{Minimal models and formality} \label{definitions}

We recall some definitions and results about minimal models
\cite{H,DGMS}. Let $(A,d )$ be a {\it differential algebra}, that
is, $A$ is a graded commutative algebra over the real numbers,
with a differential $d $ which is a derivation, i.e. $d (a\cdot b)
= (d  a)\cdot b +(-1)^{\deg (a)} a\cdot (d  b)$, where $\deg(a)$
is the degree of $a$.
A differential algebra $(A,d )$ is said to be {\it minimal\/} if:
\begin{enumerate}
 \item $A$ is free as an algebra, that is, $A$ is the free
 algebra $\bigwedge V$ over a graded vector space $V=\oplus V^i$, and
 \item there exists a collection of generators $\{ a_\tau,
 \tau\in I\}$, for some well ordered index set $I$, such that
 $\deg(a_\mu)\leq \deg(a_\tau)$ if $\mu < \tau$ and each $d
 a_\tau$ is expressed in terms of preceding $a_\mu$ ($\mu<\tau$).
 This implies that $d  a_\tau$ does not have a linear part, i.e., it
 lives in ${\bigwedge V}^{>0} \cdot {\bigwedge V}^{>0} \subset {\bigwedge V}$.
\end{enumerate}

We shall say that a minimal differential algebra $(\bigwedge V,d
)$ is a {\it minimal model} for a connected differentiable
manifold $M$ if there exists a morphism of differential graded
algebras $\rho\colon {(\bigwedge V,d )}\longrightarrow {(\Omega
M,d )}$, where $\Omega M$ is the de Rham complex of differential
forms on $M$, inducing an isomorphism
$\rho^*\colon H^*(\bigwedge V)\longrightarrow H^*(\Omega M,d )= H^*(M)$
on cohomology.

If $M$ is a simply connected manifold (or, more generally, if $M$
is a nilpotent space, i.e., $\pi_1(M)$ is nilpotent and it acts
nilpotently on $\pi_i(M)$ for $i\geq2$), then the dual of the real
homotopy vector space $\pi_i(M)\otimes \RR$ is
isomorphic to $V^i$ for any $i$. Halperin in~\cite{H} proved that
any connected manifold $M$ has a minimal model unique up to
isomorphism, regardless of its fundamental group.

A minimal model $(\bigwedge V,d )$ of a manifold $M$ is said to be
{\it formal}, and $M$ is said to be {\it formal}, if there is a
morphism of differential algebras $\psi\colon {(\bigwedge V,d
)}\longrightarrow (H^*(M),d =0)$ that induces the identity on
cohomology. An alternative way to look at this is the following:
the above property means that $(\bigwedge V,d )$ is a minimal
model of the differential algebra $(H^*(M),0)$. Therefore $(\Omega
M,d )$ and $(H^*(M),0)$ share their minimal model, i.e., one can
obtain the minimal model of $M$ out of its real cohomology
algebra. When $M$ is nilpotent, the minimal model encodes its real
homotopy type, so formality for $M$ is equivalent to saying that
its real homotopy type is determined by its real cohomology
algebra.

In order to detect non-formality, we have Massey products. Let us
recall its definition. Let $M$ be a (not necessarily simply
connected) manifold and let $a_i \in H^{p_i}(M)$, $1 \leq i\leq
3$, be three cohomology classes such that $a_1\cup a_2=0$ and
$a_2\cup a_3=0$. Take forms $\alpha_i$ in $M$ with
$a_i=[\alpha_i]$ and write $\alpha_1\wedge \alpha_2= d \xi$,
$\alpha_2\wedge \alpha_3=d \eta$. The Massey product of the
classes $a_i$ is defined as
  $$
  \langle a_1,a_2,a_3 \rangle  = [ \alpha_1 \wedge \eta+(-1)^{p_1+1} \xi
  \wedge \alpha_3] \in \frac{H^{p_1+p_2+p_3-1}(M)}{a_1
  \cup H^{p_2+p_3-1}(M) + H^{p_1+p_2-1}(M)\cup a_3} \; .
  $$

We have the following result, for whose proof we refer
to~\cite{DGMS,Ta}.

\begin{lemma} \label{lem:criterio1}
 If $M$ has a non-trivial Massey product then $M$ is non-formal.
  \hfill \QED
\end{lemma}

Therefore the existence of a non-zero Massey product is an
obstruction to the formality.

In order to prove formality, we extract the following notion from
\cite{FM}.

\begin{definition}\label{def:s-formal}
 Let $(\bigwedge V,d )$ be a minimal model of a differentiable manifold
$M$. We say that $(\bigwedge V,d )$ is \emph{$s$-formal}, or $M$
is a \emph{$s$-formal manifold} ($s\geq 0$) if for each $i\leq s$
one can get a space of generators $V^i$ of elements of degree $i$
that decomposes as a direct sum $V^i=C^i\oplus N^i$, where the
spaces $C^i$ and $N^i$ satisfy the three following conditions:
 \begin{enumerate}
 \item $d (C^i) = 0$,
 \item the differential map $d \colon N^i\longrightarrow \bigwedge V$ is
   injective,
 \item any closed element in the ideal $I_s=
   I(N^{\leq s})$, generated by
   $N^{\leq s}$ in $\bigwedge V^{\leq s}$, is exact in $\bigwedge V$.
 \end{enumerate}
\end{definition}

The condition of $s$-formality is weaker than that of formality.
However we have the following positive result proved in \cite{FM}.

\begin{theorem}\label{thm:criterio2}
Let $M$ be a connected and orientable compact differentiable
manifold of dimension $2n$ or $(2n-1)$. Then $M$ is formal if and
only if is $(n-1)$-formal. \hfill \QED
\end{theorem}

This result is very useful because it allows us to check that a
manifold $M$ is formal by looking at its $s$-stage minimal model,
that is, $\bigwedge V^{\leq s}$. In general, when computing the
minimal model of $M$, after we pass the middle dimension, the
number of generators starts to grow quite dramatically. This is
due to the fact that Poincar\'e duality imposes that the Betti
numbers do not grow and therefore there are a large number of cup
products in cohomology vanishing, which must be killed in the
minimal model by introducing elements in $N^i$, for $i$ above the
middle dimension. This makes Theorem \ref{thm:criterio2} a very
useful tool for checking formality in practice. For instance, we
have the following results, whose proofs we include for
completeness.

\begin{theorem}[\cite{Mi,FM}] \label{thm:FM}
Let $M$ be a $(k-1)$-connected compact orientable manifold of
dimension less than or equal to $(4k-2)$, $k\geq 1$. Then $M$ is
formal.
\end{theorem}

\begin{proof}
Since $M$ is $(k-1)$-connected, a minimal model $(\bigwedge V,d)$
of $M$ must satisfy $V^i = 0$ for $i\leq {k-1}$ and $V^k = C^k$
(i.e.,\ $N^k=0$). Therefore the first non-zero differential, being
decomposable, must be $d: V^{2k-1} \to V^k \cdot V^k$. This
implies that $V^j = C^j$ (i.e.,\ $N^j=0$) for $k\leq j\leq
{(2k-2)}$. Hence $M$ is $2(k-1)$-formal. Now, using
Theorem~\ref{thm:criterio2} we have that $M$ is formal.
\end{proof}

\begin{theorem}[\cite{Gil}] \label{thm:Gil}
Let $M$ be a $(k-1)$-connected compact manifold of dimension less
than or equal to $4k$, $k>1$, with $b_k=1$. Then $M$ is formal.
\end{theorem}

\begin{proof}
The minimal model $(\bigwedge V,d)$ of $M$ has $V^{<k}=0$ and
$V^k=C^k=\la \xi\ra$ one-dimensional, generated by a closed element $\xi$.
The first non-zero differential is $d: V^{2k-1} \to V^k \cdot
V^k$. This implies that $N^j=0$ for $k\leq j\leq {(2k-2)}$.
Thus $M$ is $(2k-2)$-formal. Now if $k$ is odd, then $\xi\cdot\xi=0$, so $N^{2k-1}=0$. Hence $M$ is
$(2k-1)$-formal and, since $n\leq 4k$, Theorem~\ref{thm:criterio2}
gives the formality of  $M$.

If $k$ is even, then either $N^{2k-1}=0$ and $M$ is formal as
above, or $N^{2k-1}=\la\eta\ra$ with $d\eta=\xi^2$. In this case
if $z\in I(N^{2k-1})$ is closed, write $z=\eta z_1$, $z_1 \in
\bigwedge (C^{\leq (2k-1)})$, and $0=dz = \xi^2
z_1$ which implies $z_1=0$ and hence $z=0$. Therefore $M$ is
$(2k-1)$-formal and, by Theorem~\ref{thm:criterio2} again, formal.
\end{proof}

Regarding the Massey products, we have the following refinement of
Lemma \ref{lem:criterio1}, which follows from \cite[Lemma
2.9]{FM}.

\begin{lemma}\label{lem:Massey-sformal}
Let $M$ be an  $s$-formal manifold. Suppose that there are three
cohomology classes $\alpha_i\in H^{p_i}(M)$, $1\leq i\leq 3$, such
that the Massey product $\la \alpha_1,\alpha_2,\alpha_3\ra$ is
defined. If $p_1+p_2\leq {s+1}$ and $p_{2}+ p_3\leq {s+1}$, then
$\la \alpha_1,\alpha_2,\alpha_3\ra$ vanishes.  \hfill \QED
\end{lemma}

\section{Non-formal examples with $b_k=1$} \label{construction1}

Let $k\geq 1$. In this section, we are going to give a
construction of a $(k-1)$-connected compact orientable non-formal
manifold with $b_k=1$ of any dimension $n\geq 4k+1$. The examples
that we are going to construct follow the pattern of those in
\cite{Ru} (see also \cite[Example 5]{FM}).

Consider a wedge of two spheres $S^{k}\vee S^{k+1} \subset
{\RR}^{k+2}$. Let $a\in \pi_k(S^{k}\vee S^{k+1})$ be the image of
the generator of $\pi_{k}(S^{k})$ by the inclusion of $S^{k}\inc
S^{k}\vee S^{k+1}$ and let $b\in \pi_{k+1}(S^{k}\vee S^{k+1})$ be
the image of the generator of $\pi_{k+1}(S^{k+1})$ by the
inclusion of $S^{k+1}\inc S^{k}\vee S^{k+1}$. The iterated
Whitehead product $\gamma =[a,[a,b]]\in \pi_{3k-1}(S^{k}\vee
S^{k+1})$ yields a map $\gamma: S^{3k-1} \too S^{k}\vee S^{k+1}$.
Let
  $$
  C_k= (S^{k}\vee S^{k+1}) \cup_\gamma D^{3k}
  $$
be the mapping cone of $\gamma$, where $D^{3k}$ is the
$3k$-dimensional disk.  Clearly $C_k$ is $(k-1)$-connected and its
homology groups are
  $$
  \left\{ \begin{array}{llll}
  H_i(C_k)  &= 0, & \qquad & 1\leq i\leq k-1, \\
  H_k(C_k)  &= \la a \ra,  \\
  H_{k+1}(C_k)  &= \la b \ra,\\
  H_i(C_k)  &= 0, & \qquad & k+2\leq i\leq 3k-1, \\
  H_{3k}(C_k)  &= \la c\ra,
  \end{array} \right.
  $$
where $a,b$ denote the images in $H_*(C_k)$ of the respective
elements $a,b$ under the map $\pi_*(S^k \vee S^{k+1}) \to
\pi_*(C_k) \to H_*(C_k)$, where the last map is the Hurewicz
homomorphism, and $c$ is the element represented by the attached
$3k$-cell to form $C_k$ from the wedge of the two spheres.

Let $\tilde{a}\in H^k(C_k)$, $\tilde{b}\in H^{k+1}(C_k)$ and
$\tilde{c}\in H^{3k}(C_k)$ be the cohomology classes dual to
$a,b,c$, respectively. Let us describe the minimal model of $C_k$.

\begin{itemize}
 \item Suppose $k\geq 3$ and odd. By degree reasons, all cup
 products in $H^*(C_k)$ are zero. Note that $\pi_{3k-1}(S^k\vee
 S^{k+1})$ is one-dimensional and generated by $\gamma=[a,[a,b]]$.
 Therefore $\pi_{3k-1}(C_k)=0$, since $\gamma$ contracts in $C_k$. The
 minimal model of $C_k$ is thus
   \begin{eqnarray*}
   &&(\bigwedge (\la \alpha,\beta, \eta,\xi \ra \oplus V^{\geq 3k}
   ), d), \qquad \qquad \qquad \qquad \qquad \\
   && |\alpha|=k, \ |\beta|=k+1, \ |\eta|=2k, \ |\xi|=2k+1, \\
   &&  d\alpha=d\beta=0, \ d\eta= \alpha\, \beta, \ d\xi =\beta^2,
   \end{eqnarray*}
 where $V^{\geq 3k}$ is the space generated by the generators of
 degree bigger than or equal to $3k$, $\ta=[\alpha]$,
 $\tb=[\beta]$. Note that $\alpha\, \eta$ is a closed non-exact
 element, i.e., we can write $\tc=[\alpha\, \eta]$.
The triple Massey product
  $$
  \la \ta,\ta,\tb \ra= \tc\neq 0
  \in \frac {H^{3k}(C_k)}{[\alpha] \cup H^{2k}(C_k) +
 H^{2k-1}(C_k) \cup [\beta]} = H^{3k}(C_k).
   $$
Therefore by Lemma \ref{lem:criterio1}, $C_k$ is non-formal.

 \item Suppose $k$ is even. In particular $k\geq 2$ and hence all cup
 products in $H^*(C_k)$ are zero. Again $\pi_{3k-1}(S^k\vee
 S^{k+1})$ is one-dimensional and generated by $\gamma=[a,[a,b]]$,
 hence $\pi_{3k-1}(C_k)=0$. The minimal model of $C_k$ is thus
   \begin{eqnarray*}
   &&(\bigwedge (\la \alpha,\beta, \eta,\xi \ra \oplus V^{\geq 3k}
   ),d)
   , \qquad \qquad \qquad \qquad \qquad
   \\ && |\alpha|=k, \ |\beta|=k+1, \ |\eta|=2k-1, \ |\xi|=2k, \\
   && d\alpha=d\beta=0, \ d\eta= \alpha^2, \ d\xi =\alpha\,\beta,
  \end{eqnarray*}
 where $\ta=[\alpha]$,
 $\tb=[\beta]$. Now $\alpha\, \xi +\beta\,\eta$ is a closed non-exact
 element of degree $3k$. So we can write $\tc=[\alpha\, \xi +\beta\,\eta]$.
 The triple Massey product
 $$
 \la \ta,\ta,\tb \ra= \tc\neq 0
 \in \frac {H^{3k}(C_k)}{[\alpha] \cup H^{2k}(C_k) +
 H^{2k-1}(C_k) \cup [\beta]} = H^{3k}(C_k),
 $$
and by Lemma \ref{lem:criterio1}, $C_k$ is non-formal.

\item The case $k=1$ is slightly different. Here $\pi_2(S^1\vee
S^2)$ is infinitely generated by $b$, $[a,b]$, $[a,[a,b]]$,
$[a,[a,[a,b]]], \ldots$ Therefore $\pi_2(C_1)$ is $2$-dimensional
generated by $b, [a,b]$. So the minimal model of $C_1$ is
  \begin{eqnarray*}
   &&( \bigwedge (\la \alpha,\beta, \eta  \ra \oplus V^{\geq 3} ),
   d), \qquad \qquad \qquad \qquad \qquad \\
   && |\alpha|=1, \ |\beta|=2, \ |\eta|=2, \\
   && d\alpha=d\beta=0, \ d\eta= \alpha\,\beta.
  \end{eqnarray*}
Here $\tc=[\alpha\ \eta]$  and we have the
non-vanishing triple Massey product
 $$
 \la \ta,\ta,\tb \ra= \tc\neq 0
 \in \frac {H^{3}(C_1)}{[\alpha] \cup H^{2}(C_1) +
 H^{1}(C_1) \cup [\beta]} = H^{3}(C_1),
 $$
 proving non-formality of $C_1$.
\end{itemize}

Now we aim to construct a differentiable compact manifold using
$C_k$. For this, we use \cite[Corollary 2]{Ru} to obtain a PL
embedding $C_k\subset {\RR}^{3k+(k+2)} = {\RR}^{4k+2}$. Let $W_k$
be a closed regular neighborhood of $C_k$ in $\RR^{4k+2}$ and let
$Z_k= \bd W_k$ be its boundary. We can arrange easily that $Z_k$
is a {\em smooth\/} manifold of dimension $4k+1$.

\begin{theorem} \label{thm:mainexample}
 Suppose $k\geq 2$. Then $Z_k$ is a $(k-1)$-connected compact (orientable)
 non-formal $(4k+1)$-dimensional manifold
 with $b_k(Z_k)=1$.
\end{theorem}

\begin{proof}
Suppose $k\geq 1$ by now.
The first observation is that $\pi_i(Z_k)\iso \pi_i(W_k-C_k)\iso
\pi_i(W_k) \iso \pi_i(C_k)=0$, for $1\leq i\leq k-1$, where the
first and last isomorphism are by retraction deformation, and the
middle one because $C_k$ has codimension $k+2$ in $W_k$. Therefore
$Z_k$ is $(k-1)$-connected. Moreover $\pi_k(Z_k) \iso \pi_k(C_k)\iso
\ZZ$, hence $b_k(Z_k)=1$.

Now let us see that $Z_k$ is non-formal. For this, we need to
compute its cohomology. There is a long exact sequence
  $$
   \cdots \too H^{i}(W_k,Z_k) \stackrel{j^*}{\too} H^i(W_k)
   \stackrel{i^*}{\too} H^i(Z_k) \stackrel{\bd^*}{\too} H^{i+1}(W_k,Z_k) \too
   \cdots,
  $$
where $i:Z_k\to W_k$ and $j:W_k \to (W_k,Z_k)$ are the inclusions.
Using that $H^*(C_k)\iso H^*(W_k)$ and $H^*(W_k,Z_k)\iso
H_{4k+2-*}(W_k) \iso H_{4k+2-*}(C_k)$, the first isomorphism by
Poincar{\'e} duality, we rewrite the above sequence as
  $$
   \cdots \too H_{4k+2-i}(C_k) \stackrel{j^*}{\too} H^i(C_k)
   \stackrel{i^*}{\too} H^i(Z_k) \stackrel{\bd^*}{\too}
   H_{4k+1-i}(C_k) \too \cdots
 $$
{}From this it follows easily that $i^*$ is always injective. The
only non-trivial case is $k=1$ where there is a map $H_3(C_1)\to
H^3(C_1)$, but this is an antisymmetric map between rank one
spaces, hence the zero map.

We deduce the following cohomology groups for any $k\geq 1$:
  $$
  \left\{\begin{array}{llll}
  H^i(Z_k)  &= 0, & \qquad & 1\leq i\leq k-1, \\
  H^k(Z_k)  &= \la \check{a} \ra,  \\
  H^{k+1}(Z_k)  &= \la \check{b},\hc \ra,\\
  H^i(Z_k)  &= 0, & \qquad & k+2\leq i\leq 3k-1, \\
  H^{3k}(Z_k)  &= \la \check{c},\hb \ra, \\
  H^{3k+1}(Z_k)  &= \la \ha \ra,\\
  H^i(Z_k)  &= 0, & \qquad & 3k+2 \leq i\leq 4k, \\
  H^{4k+1} (Z_k) & =\la [Z_k]\ra,
  \end{array} \right.
  $$
where $\check{a},\check{b},\check{c}$ denote the images of
$\ta,\tb,\tc \in H^*(C_k)$ under $i^*$, and $\ha,\hb,\hc$ denote
the preimages of $\ta,\tb,\tc \in H^*(C_k)$ under $\bd^*$, and
$[Z_k]$ is the fundamental class of $Z_k$.

Now suppose that $k\geq 2$.
Let us see that $Z_k$ has a non-vanishing Massey product. As
$i^*:H^*(C_k)\to H^*(Z_k)$ is
injective, there is a Massey product $\la i^*\ta,i^*\ta,i^*\tb \ra
=i^*\tc$ and this is non-zero in
 $$
 \frac {H^{3k}(Z_k)}{\check{a}
 \cup H^{2k}(Z_k)+H^{2k-1}(Z_k) \cup \check{b}} = H^{3k}(Z_k),
 $$
using that $H^{2k}(Z_k)=0$ and $H^{2k-1}(Z_k)=0$ for $k \geq 2$
(this is the only place where the assumption $k\geq 2$ is used).
\end{proof}

For constructing higher dimensional examples, we may use lemma
\ref{lem:higher} below, but a more direct way is available, as follows.

\begin{theorem} \label{thm:mainexample2}
 Let $k\geq 1$. There are $(k-1)$-connected compact non-formal $n$-dimensional
 manifold $Z_{k,n}$ with $b_k(Z_{k,n})=1$, for any $n\geq 4k+1$, with $n\not=5k$.
\end{theorem}

\begin{proof}
 For any $n\geq 4k+1$,
 embed $C_k\subset \RR^{4k+2}\subset \RR^{n+1}$. Take a tubular
 neighborhood $W_{k,n}$ of $C_k$ in $\RR^{n+1}$
 and let $Z_{k,n}=\bd W_{k,n}$ be its boundary. The same
 argument as in the proof of Theorem \ref{thm:mainexample} proves
 that $Z_{k,n}$ is
 a $(k-1)$-connected compact orientable manifold with
 $b_k(Z_{k,n})=1$. Let us see that it is non-formal by checking
 that it has a non-vanishing triple Massey product.

Let us see that the map $i^*:H^*(C_k)\iso H^*(W_{k,n})\to
H^*(Z_{k,n})$ is injective. We have a commutative diagram
 \begin{eqnarray*}
  H^*(W_{k,n}) &\stackrel{i^*}{\too} &   H^*(Z_{k,n}) \\
   \downarrow\cong  & & \downarrow \\
  H^*(W_{k}) &\stackrel{i^*}{\too} &   H^*(Z_{k}).
  \end{eqnarray*}
Since the bottom row is injective by the proof of
Theorem \ref{thm:mainexample}, the top one is also. Thus
there is a Massey product $\la i^*\ta,i^*\ta,i^*\tb \ra =i^*\tc$.

Now assume $n\neq 5k$, and let us prove that  $i^*\ta
\cup H^{2k}(Z_{k,n})+H^{2k-1}(Z_{k,n}) \cup i^*\tb=0 \subset
H^{3k}(Z_{k,n})$. As in the proof of Theorem
\ref{thm:mainexample}, we have an exact sequence
  $$
   \cdots \too H_{n+1-i}(C_k) \stackrel{j^*}{\too} H^i(C_k)
   \stackrel{i^*}{\too} H^i(Z_{k,n}) \stackrel{\bd^*}{\too}
   H_{n-i}(C_k) \too \cdots
 $$
Then $H^{2k}(Z_{k,n})=i^*H^{2k}(C_k)$ since $H_{n-2k}(C_k)=0$.
Then $i^*\ta \cup H^{2k}(Z_{k,n})=0$. Analogously
$H^{2k-1}(Z_{k,n}) \cup i^*\tb =0$, unless $n=5k-1$.
If $n=5k-1$ (in particular, $k\geq 2$) then
$H^{2k-1}(Z_{k,n})=({\bd^*})^{-1}H_{3k}(C_k)$,
generated by an element $\hat{c}$.
Now $i^*\tb \cup \hat{c}=0\in H^{3k}(Z_{k,n})$, by Poincar\'e
duality since $i^*\tb \cup \hat{c}\cup \hat{c}=0$ (as $\hat{c}$
has odd degree). So $H^{2k-1}(Z_{k,n}) \cup i^*\tb =0$.

It follows that $i^*\ta \cup H^{2k}(Z_{k,n})+H^{2k-1}(Z_{k,n}) \cup i^*\tb=0$
and so the above Massey product is non-zero in
 $$
 \frac {H^{3k}(Z_{k,n})}{i^*\ta \cup
 H^{2k}(Z_{k,n})+H^{2k-1}(Z_{k,n}) \cup i^*\tb} = H^{3k}(Z_{k,n}).
 $$
\end{proof}

\section{Nilpotency of the constructed examples}
\label{nilpotency}

In section \ref{construction1} we provide with some examples of non-simply
connected manifolds, the manifolds $Z_{1,n}$, $n\geq 5$.
As we are studying a rational homotopy property,
it is a natural question whether or not they are nilpotent spaces.
Here we collect here a rather non-conclusive collection of remarks
on this question.

The fundamental group of $Z_{1,n}$ is $\pi_1(Z_{1,n}) =\la a\ra \iso \ZZ$,
abelian. To check nilpotency of $Z_{1,n}$, need to describe the action of $a$
on the higher homotopy groups $\pi_i(Z_{1,n})$.
For instance, suppose $n\geq 6$
(so we already know that $Z_{1,n}$ is non-formal
by Theorem \ref{thm:mainexample2}). Then $\pi_2(Z_{1,n}) \iso
\pi_2(C_1)$. Since $[a,[a,b]] =0\in \pi_2(C_1)$, the action
of $a$ on this homotopy group is nilpotent.
In general, the nilpotency of the action of $a$ on higher
homotopy groups $\pi_i(Z_{1,n})$, $i>2$ reduces to two
issues:
 \begin{itemize}
  \item The nilpotency of $C_1$. It is not clear whether the action
   of $a$ on the higher homotopy groups $\pi_i(C_1)$, $i>2$, is nilpotent,
   although this seems to be the case. Let us do the case $i=3$.

The Quillen model \cite{Ta} of $C_1$ is $(\LL(a,b,c),\bd)$,
where $\LL= \LL(a,b,c)$ is the free Lie algebra generated by
elements $a,b,c$ of degrees $0,1,2$ and with $\bd a=0$, $\bd b=0$, $\bd c=
[a,[a,b]]$. Then $\pi_i(C_1) \iso H^{i-1}(\LL(a,b,c),\bd)$.
Consider the map $p:\LL \to \LL$, $p(x)=[a,x]$ of degree $0$. This is
a derivation. Moreover, a basis for the elements of degree $2$ in $\LL$
is $\{ p^j(c), p^j([b,b]); j\geq 0\}$. Let $z= \sum \lambda_j p^j(c)
+\mu_j p^j([b,b])$ be a closed element, defining a homology class
$\bar{z} \in H^{2}(\LL,\bd)= \pi_3(C_1)$. Hence
$0 =\bd z= \sum \lambda_j p^{j+2}(c)$, so $\lambda_j=0$ for all $j\geq 0$.
So we can write $\bar{z}=\sum \mu_j \overline{p^j[b,b]}$.
The map $p$ descends to homology and $p^2 (\bar{b})= \overline{ p^2 (b)}
=\overline{\bd c}=0$. As $p$ is a derivation,
$p^3 \left( \overline{[b,b]} \right) = 0$.
Hence $p:\pi_3(C) \to \pi_3(C)$ is nilpotent.

\item Knowing the nilpotency of $C_1$, prove the nilpotency of $Z_{1,n}$.
 For $i\leq n-4$, $\pi_i(Z_{1,n})\iso \pi_i(C_1)$. So nilpotency of the
 action of $a$ on
$\pi_i(Z_{1,n})$ would follow from the corresponding statement for $C_1$.
Now suppose $i=n-3$. Then as $Z_{1,n}$ is of the homotopy type of
$W_{1,n} - C_1 \subset \RR^{n+1}$, we have that $\pi_{n-3}(Z_{1,n})$ is
generated by $\pi_{n-3}(C_1)$ and the $S^{n-3}$-fiber $f$ of
the projection $Z_{1,n}\to C_1$ over a smooth point of $C_1$.
It is not clear that the action of $a$ on $f$ is nilpotent.
Even more difficult is the case $i>n-3$.
\end{itemize}

\begin{remark}
In \cite{FM3} it is claimed that the examples constructed
therein of non-formal
manifolds with $b_1=1$ are not nilpotent (see \cite[Section 5]{FM3}).
In \cite[Lemma 9]{FM3}
it is proved that in our circumstances
the action of any non-zero element $a\in \pi_1(M)$ on $\pi_2(M)$ is
not trivial. This action is given as
$h_a:\pi_2(M)\to \pi_2(M)$, $h_a(A)=[a,A] +A=p(A)+A$.
Therefore $h^k(A) \neq A$ for any $A\neq 0$, $k\geq 1$. However this does not mean
that the action is not nilpotent, since nilpotency means that $p^N=
(h-Id)^N=0$, for some large $N$. It may happen that the examples of \cite{FM3}
are nilpotent, though the authors do not know the answer.
\end{remark}

\section{Non-formal examples with $b_k=2$} \label{construction2}

Let $k\geq 1$. In this section, we are going to give a
construction of a $(k-1)$-connected compact orientable non-formal
manifold with $b_k=2$ of any dimension $n\geq 4k$, by modifying
slightly the construction in Section \ref{construction1}.

Consider now a wedge of two spheres $S^{k}\vee S^{k} \subset
{\RR}^{k+1}$. Let $a,b\in \pi_k(S^{k}\vee S^{k})$ be the image of
the generators of $\pi_{k}(S^{k})$ by the inclusions of $S^{k}\inc
S^{k}\vee S^{k}$ as the first and second factors, respectively.
The iterated Whitehead product $\gamma =[a,[a,b]]\in
\pi_{3k-2}(S^{k}\vee S^{k})$ yields a map $\gamma: S^{3k-2} \too
S^{k}\vee S^{k}$ and a $(3k-1)$-dimensional CW-complex
  $$
  C_k'= (S^{k}\vee S^{k}) \cup_\gamma D^{3k-1}.
  $$
Clearly $C_k'$ is $(k-1)$-connected and the only non-zero homology
groups are
  $$
  \left\{ \begin{array}{llll}
  H_k(C_k')  &= \la a,b \ra,  \\
  H_{3k-1}(C_k')  &= \la c\ra,
  \end{array} \right.
  $$
where $a,b$ denote the images of the elements $a,b$ under the map
$\pi_*(S^k \vee S^{k}) \to \pi_*(C_k') \to H_*(C_k')$, and $c$ is
the element represented by the attached $(3k-1)$-cell to form
$C_k'$ from the wedge of the two spheres. Let $\ta,\tb \in
H^k(C_k')$ and $\tc \in H^{3k-1}(C_k')$ be the cohomology classes
dual to $a,b,c$, respectively. Then it is easy to see, as before,
that the triple Massey product
  $$
  \la \ta,\ta,\tb \ra= \tc\neq 0
  \in \frac {H^{3k-1}(C_k')}{\ta \cup H^{2k-1}(C_k') +
  H^{2k-1}(C_k') \cup \tb} = H^{3k-1}(C_k').
  $$
Therefore by Lemma \ref{lem:criterio1}, $C_k'$ is non-formal.

Now we take a PL embedding $C_k'\subset {\RR}^{3k-1+(k+1)} =
{\RR}^{4k}\subset \RR^{n+1}$, let $W_{k,n}'$ be a closed regular
neighborhood of $C_k'$ in $\RR^{n+1}$. Let $Z_{k,n}'= \bd
W_{k,n}'$ be its boundary. We can arrange easily that $Z_{k,n}'$
is a smooth manifold of dimension $n$.

\begin{theorem} \label{thm:mainexample3}
 $Z_{k,n}'$ is a $(k-1)$-connected compact orientable non-formal
 $n$-dimensional manifold with $b_k(Z_{k,n}')=2$, for any $n\geq 4k$.
\end{theorem}

\begin{proof}
Note that the codimension of $C_k'$ in $W_{k,n}'$ is
$n+1-(3k-1)\geq k+2$. Therefore $\pi_i(Z_{k,n}')\iso
\pi_i(W_{k,n}'-C_k')\iso \pi_i(W_{k,n}') \iso \pi_i(C_k')=0$, for
$1\leq i\leq k-1$, so $Z_{k,n}'$ is $(k-1)$-connected. Also
$\pi_k(Z_{k,n}')\iso \pi_k(C_k')$, so $b_k(Z_{k,n}')=b_k(C_k')=2$
(this is the reason for the necessity of the condition $n\geq
4k$).

Now let us see that $Z_{k,n}'$ is non-formal. For this, we need to
compute its cohomology. There is a long exact sequence
  $$
   \cdots \too H^{i}(W_{k,n}',Z_{k,n}') \stackrel{j^*}{\too}
   H^i(W_{k,n}') \stackrel{i^*}{\too} H^i(Z_{k,n}')
   \stackrel{\bd^*}{\too} H^{i+1}(W_{k,n}',Z_{k,n}') \too
   \cdots,
  $$
where $i:Z_{k,n}'\to W_{k,n}'$ and $j:W_{k,n}' \to
(W_{k,n}',Z_{k,n}')$ are the inclusions. Use that $H^*(C_{k}')\iso
H^*(W_{k,n}')$ and $H^*(W_{k,n}',Z_{k,n}')\iso H_{n+1-*}(W_{k,n}')
\iso H_{n+1-*}(C_k')$. Hence we rewrite the above sequence as
  $$
   \cdots \too H_{n+1-i}(C_k') \stackrel{j^*}{\too} H^i(C_k')
   \stackrel{i^*}{\too} H^i(Z_{k,n}') \stackrel{\bd^*}{\too}
   H_{n-i}(C_k') \too \cdots
 $$

The map $i^*:H^*(C_k')\iso H^*(W_{k,n}') \to H^*(Z_{k,n}')$
is injective. In the case $n=4k$, $H^*(W_{k,4k}') \to
H^*(Z_{k,4k}')$ is injective, because $H_{4k+1-i}(C_k')=0$
for $i=k,3k-1$. Now for $n\geq 4k$, we have a
commutative diagram
 \begin{eqnarray*}
  H^*(W_{k,n}) &\stackrel{i^*}{\too} &   H^*(Z_{k,n}) \\
   \downarrow\cong  & & \downarrow \\
  H^*(W_{k,4k}) &\stackrel{i^*}{\too} &   H^*(Z_{k,4k}),
  \end{eqnarray*}
where the bottom row is injective, hence the top one is also.
This proves the injectivity of $i^*$.
Thus $H^{k}(Z_{k,n}')=\la i^*\ta,i^*\tb \ra$ and there is a
well-defined Massey product $\la i^*\ta,i^*\ta,i^*\tb \ra
=i^*\tc$.

Finally let us see that $i^*\ta
 \cup H^{2k-1}(Z_{k,n}')+H^{2k-1}(Z_{k,n}') \cup i^*\tb=0\subset
H^{3k-1}(Z_{k,n}')$.
First suppose $n \neq 5k-2$. Then $H^{2k-1}(Z_{k,n}') = i^*H^{2k-1}(C_k')$,
since $H_{n-2k+1}(C_k') =0$. So
$i^*\ta \cup H^{2k-1}(Z_{k,n}') = i^* \ta \cup i^*H^{2k-1}(C_k')=0$.
Analogously, $H^{2k-1}(Z_{k,n}') \cup i^*\tb=0$.

The remaining case is $n=5k-2$. Then $k>1$ since $n\geq 4k$.
So $H^{2k-1}(Z_{k,n}') =
({\bd^*})^{-1} H_{3k-1}(C_k')$. This is generated by an element
$\hat{c}$. Now $i^* a \cup \hat{c}=0\in H^{3k-1}(Z_{k,n}')$,
using Poincar\'e duality and
$i^* a \cup \hat{c}\cup\hat{c}=0$ (the degree of $\hat{c}$ is odd).
So $i^*\ta \cup H^{2k-1}(Z_{k,n}') =0$. Analogously,
$H^{2k-1}(Z_{k,n}') \cup \tb=0$.

This implies that the Massey product
$\la i^*\ta,i^*\ta,i^*\tb \ra =i^*\tc$ is
non-zero in
 $$
 \frac {H^{3k-1}(Z_{k,n}')}{i^*\ta
 \cup H^{2k-1}(Z_{k,n}')+H^{2k-1}(Z_{k,n}') \cup i^*\tb}
 = H^{3k-1}(Z_{k,n}').
 $$
\end{proof}

\section{Proof of Theorem \ref{thm:main}} \label{proof}

We start with some elementary lemmata.

\begin{lemma} \label{lem:connected}
 If $M$ has a non-trivial Massey product and $N$ is any smooth
 manifold, then $M\# N$ has a non-trivial Massey product.
\end{lemma}

\begin{proof}
  Let $\langle a_1,a_2,a_3 \rangle$ be a non-zero Massey product on $M$,
  $a_i \in H^{p_i}(M)$, $1 \leq i\leq 3$. Since $p_i>0$, it is easy
  to arrange that $a_i=[\alpha_i]$, $1\leq i\leq 3$,
  $\alpha_1\wedge \alpha_2= d \xi$ and $\alpha_2\wedge \alpha_3=d \eta$
  where $\alpha_i$, $1\leq i\leq 3$, $\xi$ and $\eta$ are forms vanishing
  on a given disc in $M$ (see \cite{BoTu}). Using this disk for
  performing the connected sum, we see that we can define the forms
  $\alpha_i$, $1\leq i\leq 3$, $\xi$ and $\eta$ in $M\# N$ by extending by
  zero. Let $a_i'=[\alpha_i] \in H^{p_i}(M\# N)$ be the cohomology classes
  thus defined. It follows easily that
  $$
  \langle a_1',a_2',a_3' \rangle  = \langle a_1,a_2,a_3 \rangle
  \in \frac{H^{p_1+p_2+p_3-1}(M)}{a_1
  \cup H^{p_2+p_3-1}(M) + H^{p_1+p_2-1}(M)\cup a_3}
  \subset \qquad \qquad
  $$
  $$
  \qquad\qquad\qquad\qquad \subset \frac{H^{p_1+p_2+p_3-1}(M\# N)}{a_1'
  \cup H^{p_2+p_3-1}(M\# N) + H^{p_1+p_2-1}(M\# N)\cup a_3'}
  $$
  is non-zero.
\end{proof}

\begin{lemma} \label{lem:higher}
 Let $M$ be an $n$-dimensional manifold with a non-trivial
 Massey product $\langle a_1,a_2,a_3 \rangle$,
 $a_i \in H^{p_i}(M)$, $1 \leq i\leq 3$, with $p_1+p_2+p_3<n$.
 Then the $(n+1)$-dimensional manifold $N=(M \times
 S^1) \#_{S^1} \, S^{n+1}$ has a non-trivial Massey product.
 Moreover, if $M$ is $(k-1)$-connected then so is $N$ and $b_k(N)=b_k(M)$.
\end{lemma}

\begin{proof}
  Note that
  $$
  N=(M \times S^1) \#_{S^1} \, S^{n+1} =((M - D^n)\times S^1) \cup_{S^{n-1}
  \times S^1} (S^{n-1} \times D^2),
  $$
  as there is only one way to embed $S^1$ in $S^{n+1}$ since $n\geq 3$;
  otherwise $M$ cannot have non-trivial Massey products.
  Any cohomology class $a\in H^*(M)$ of positive degree has a representative
  vanishing on the disc $D^n$. Therefore it defines in a natural way a cohomology class
  on $N$, giving a map $H^*(M) \to H^*(N)$.
  A Mayer-Vietoris argument gives that the cohomology of $N$ is
  $H^k(N)=H^k(M)$, $k=0,1$, $H^k(N)=
  H^k(M) \oplus H^{k-1}(M)\cdot [\eta]$, $2\leq k\leq n-1$, $H^k(N)=
  H^{k-1}(M)\cdot [\eta]$, $k=n,n+1$ (where $[\eta]$ is the generator of $H^1(S^1)$).
  {}From this it follows easily the last sentence of the statement.

  Now write $a_i=[\alpha_i]$, $1 \leq i\leq 3$, with
  $\alpha_1\wedge \alpha_2= d \xi$ and $\alpha_2\wedge \alpha_3=d \mu$.
  We arrange that $\alpha_i$, $1\leq i\leq 3$, $\xi$ and $\mu$ are
  forms on $M$ which vanish on the given disc $D^n \subset M$.
  This yields
  a Massey product $\langle a_1',a_2',a_3' \rangle$ on $N$,
  $a_i'=[\alpha_i] \in H^{p_i}(N)$.
  Since the map $H^*(M) \to H^*(N)$ is injective for $*<n$, it follows that
  this Massey product is non-zero in $N$.
\end{proof}

\bigskip

\noindent {\bf Proof of Theorem \ref{thm:main}\ :} \
First let us address the \textit{only if}\/ part of the theorem.
Let $M$ be a compact orientable $(k-1)$-connected $n$-manifold. If
$n\leq 4k-2$ then $M$ is formal by Theorem \ref{thm:FM}. If $M$
has $b_k=1$ and its dimension is $n\leq 4k$ then $M$ is formal by
Theorem \ref{thm:Gil}. If $M$ has $b_k=0$ (this means that either
$M$ is $k$-connected, or else that $\pi_k(M)$ is torsion) then the
minimal model of $M$ is of the form $(\bigwedge V^{\geq
(k+1)},d)$. The argument of the proof of Theorem \ref{thm:FM}
proves the formality of $M$ if $n\leq 4k+2$.

For the \textit{if}\/ part of the theorem, we have to give
constructions of non-formal compact orientable $(k-1)$-connected
$n$-manifolds for any $b_k=b\geq 0$ and $n\geq \max \{ 4k-1 ,
4k+3-2b\}$.

\begin{itemize}
 \item Case $b_k=0$ and $n\geq 4k+3$. We need examples of
 $(k-1)$-connected non-formal $n$-manifolds with $b_k=0$.
 For instance, take the $k$-connected non-formal
 $n$-manifolds provided by \cite{Ru},
 since $n\geq 4(k+1)-1$.

 \item Case $b_k=1$ and $n\geq 4k+1$. The manifold $Z_{k,n}$ provided
 by Theorem \ref{thm:mainexample2} covers this case when $n\not= 5k$.
 For $k=1$, in \cite{FM3} are given examples of non-formal
 $5$-dimensional manifolds with first Betti number $b_1=1$.
 For $n=5k$ with $k \geq 2$, it is sufficient to consider the manifold
 $(Z_{k,5k-1}\times S^1) \#_{S^1} \, S^{5k}$, by
 Lemma \ref{lem:higher}.

 \item Case $b_k=2$ and $n= 4k-1$. For the specific case $k=2$, where
we have $n=7$ and $b_2=2$, Oprea \cite{O} constructs an example of
a compact non-formal manifold as the total space of a $S^3$-bundle
over $S^2\times S^2$ with Euler class $1$ (actually Oprea gives a
different construction, but it is easily seen to reduce to the
above description). This was the first example of a non-formal
simply connected compact manifold of dimension $7$.

This construction is generalized by Cavalcanti \cite[Example
1]{Gil} for any $k\geq 1$. The total space of a $S^{2k-1}$-bundle
over $S^k\times S^k$ with Euler class $1$, is a non-formal
$(k-1)$-connected compact orientable $(4k-1)$-dimensional manifold
with $b_k=2$.

 \item Case $b_k=2$ and $n\geq 4k$. The manifold $Z_{k,n}'$ provided
 by Theorem \ref{thm:mainexample3} covers this case.
 Otherwise, apply Lemma \ref{lem:higher} repeatedly to the previous
 example.

\item Case $b_k > 2$ and $n\geq 4k-1$. Let $Z$ be a non-formal
$(k-1)$-connected orientable compact $n$-manifold with $b_k=2$.
Consider $Z \# (b_{k}-2) \, (S^{k+1}\times S^{n-k-1})$, which is
non-formal by Lemma \ref{lem:connected}.

\end{itemize}

\hfill \QED

\begin{remark} \label{rem:Gil}
Cavalcanti \cite{Gil} gives also examples of non-formal
$(k-1)$-connected compact orientable $n$-dimensional manifolds
with $b_k=1$, for $n\geq 4k+1$. The total space of a
$S^{2k+2i-1}$-bundle over $S^k\times S^{k+2i}$ with Euler class
$1$ is non-formal, of dimension $4k+4i-1$ with $b_k=1$, for $i>0$.
This covers the case $b_k=1$, $n=4k+4i-1\geq 4k+3$ in the list above
(in \cite[Example 1]{Gil} it is shown an improvement to cover also
the case $b_k=1$, $n=4k+4i$, with $i>0$). Note that this method does
not gives examples for the minimum possible value $n=4k+1$.
\end{remark}

\begin{remark} \label{rem:sformality}
The examples $Z_{k,n}$ are $(k-1)$-connected with $b_k=1$. Hence
by the proof of Theorem \ref{thm:FM}, $Z_{k,n}$ is
$(2k-1)$-formal. It is not $2k$-formal by Lemma
\ref{lem:Massey-sformal}. Note that $n=4k+1$ is the smallest
dimension in which this can happen by Theorem \ref{thm:criterio2}.

The examples $Z_{k,n}'$ are $(k-1)$-connected. Hence by the proof
of Theorem \ref{thm:Gil}, $Z_{k,n}'$ is $(2k-2)$-formal. It is not
$(2k-1)$-formal by Lemma \ref{lem:Massey-sformal}. Again Theorem
\ref{thm:criterio2} says that $n=4k-1$ is the smallest dimension
in which this can happen.
\end{remark}

\section{Non-formal manifolds with small Betti numbers $b_k$
and $b_{k+1}$} \label{surgery}

A natural question that arises from the proof of Theorem
\ref{thm:mainexample} is whether there are examples of compact
non-formal $k$-connected $n$-manifolds, $n\geq 4k+1$ with $b_k=1$
and $b_{k+1}$ as small as possible.
Our examples with $n>4k+1$ satisfy that $b_{k+1}=1$. But the
examples with $n=4k+1$ satisfy that $b_{k+1}=2$.

\begin{lemma}
 If $M$ is a $(k-1)$-connected compact orientable $n$-manifold
 with $b_k=1$, $b_{k+1}=0$ and $n\leq 4k+2$, $k\geq 1$, then $M$ is formal.
\end{lemma}

\begin{proof}
 Work as in the proof of Theorem \ref{thm:Gil} to conclude that $M$ is
 $(2k-1)$-formal. For $k>2$, we have $V^{k+1}=0$, since $b_{k+1}=0$. So there is
 no product to be killed in $V^k \cdot V^{k+1}$ and hence $N^{2k}=0$
 (but maybe $C^{2k}\not=0$) proving  that $M$ is $2k$-formal.
 If $k=2$, then $V^2=C^2=\la \xi \ra$. If moreover $V^3=0$ then we
 conclude the $4$-formality as above. Otherwise $V^3=N^3=\la \eta\ra$ with
 $d\eta=\xi^2$. So there is nothing closed in $V^2 \cdot V^3$, and we have  $N^4=0$.
Therefore, $M$ is $4$-formal. Finally,
 if $k=1$, then $V^1=\la \xi \ra$ and $V^2=0$, hence $M$ is $2$-formal.
 Now use Theorem \ref{thm:criterio2} to get the formality of $M$.
\end{proof}

\bigskip

We end up with the following

\noindent\textbf{Question:}
 Let $M$ be a $(k-1)$-connected compact orientable $n$-manifold
 with $b_k=1$, $b_{k+1}=1$ and $n\leq 4k+2$, $k \geq 1$. Is it $M$ is formal?

\bigskip

\noindent {\bf Acknowledgments.}
Thanks are  due to the Organizing Committee who worked
so hard to make the conference in Belgrade a success.
We are grateful to Gil Cavalcanti for pointing us to
a version of Lemma \ref{lem:higher}.
Second author is grateful
to the Departamento de Matem\'atica Pura da Universidade do Porto, where
part of this work was done.
This work has been partially
supported through grant
grants MCyT (Spain) Project MTM2004-09070-C03-01 and
UPV 00127.310-E-15909/2004 .

\vspace{0.15cm}

{\small \noindent{\sf M. Fern\'andez:} Departamento de
Matem\'aticas, Facultad de Ciencia y Tecnolog\'{\i}a, Universidad
del Pa\'{\i}s Vasco, Apartado 644, 48080 Bilbao, Spain.

{\sl E-mail:} marisa.fernandez@ehu.es

\vspace{0.15cm}

\noindent{\sf V. Mu\~noz:} Departamento de Matem\'aticas, Consejo
Superior de Investigaciones Cient{\'\i}ficas, C/ Serrano 113bis, 28006
Madrid, Spain.

{\sl E-mail:} vicente.munoz@imaff.cfmac.csic.es

}

\end{document}